\documentclass[10pt]{amsart}
\usepackage{amsfonts,amsmath,amssymb,amsthm}
\usepackage{color}
\usepackage{graphicx}
\usepackage{hyperref}
\hypersetup{backref,colorlinks=true,citecolor=blue,linkcolor=blue}
\newtheorem{thm}{Theorem}[section]

\newtheorem{cor}[thm]{Corollary}
\newtheorem{lmm}[thm]{Lemma}
\newtheorem{prp}[thm]{Proposition}

\theoremstyle{remark}
\newtheorem*{rmk}{Remark}

\newcommand{\Mod}[1]{\mathrm{Mod}({#1})}
\newcommand{\Homeo}[1]{\mathrm{Homeo}^{+}({#1})}
\newcommand{\C}[1]{\mathcal{C}({#1})}
\newcommand{\V}[1]{\mathcal{C}_{0}({#1})}
\newcommand{\G}[1]{\mathcal{C}_{1}({#1})}

\newcommand{\grp}[1]{\langle{#1}\rangle}

\newcommand{\dis}[2]{\mathrm{d}_{{#1}}({#2})}
\newcommand{\Z}{\mathbb{Z}}
\newcommand{\R}{\mathbb{R}}
\newcommand{\nada}{\emptyset}
\newcommand{\del}{\partial}
\newcommand{\sez}{\Longrightarrow}
\newcommand{\QED}{\hfill$\Box$\medskip}

\newcommand{\goesto}{\rightarrow}

\title[Uniform uniform exponential growth of mapping class subgroups]{Uniform uniform exponential growth of subgroups of the mapping~class~group}
\author{Johanna Mangahas}
\address{University of Michigan}
\email{mangahas@umich.edu}
\thanks{The author is partially supported by NSF RTG grant \#0602191}

\begin{document}
\maketitle

\begin{abstract}
Let $\Mod{S}$ denote the mapping class group of a compact, orientable surface $S$.  We prove that finitely generated subgroups of $\Mod{S}$ which are not virtually abelian have uniform exponential growth with minimal growth rate bounded below by a constant depending only, and necessarily, on $S$.  For the proof, we find in any such subgroup explicit free group generators which are ``short'' in any word metric.  Besides bounding growth, this allows a bound on the return probability of simple random walks.
\end{abstract}

\section{Introduction}

\noindent Let $G$ be a group generated by a finite set $A$. Denote by $b(G,A,n)$ the number of elements of $G$ which may be written as words in $A$ with length less than $n$.  The exponential growth rate of $G$ with respect to $A$ is
\[h(G,A) = \limsup_{n \to \infty} {\log(b(G,A,n)) \over n}\]
and the minimal exponential growth rate of $G$ is
\[h(G) = \inf \{h(G,A) : A \text{ generates } G\}.\]
$G$ is said to have exponential growth if $h(G,A) > 0$ for some and hence all $A$, and uniform exponential growth if $h(G) > 0$.

Exponential growth does not generally imply uniform exponential growth \cite{W}, but Eskin, Mozes, and Oh proved this implication holds for linear groups over a field with characteristic zero \cite{EMO}.  Breuillard and Gelander strengthened the uniformity: for example, they proved that for a finitely generated, discrete subgroup $G < GL_n(\R)$ which is not virtually nilpotent, $h(G)$ has a lower bound depending only on $n$ \cite{BG}.  In this note we mirror this property for mapping class groups of surfaces.  It is not known whether mapping class groups are linear, except in special cases: most notably, Krammer and Bigelow proved linearity of braid groups (mapping class groups of punctured spheres) \cite{K}, \cite{Bi}.  However, many linear group results inspire mapping class group analogues.  For example, McCarthy and Ivanov independently proved that mapping class groups satisfy a version of the Tits alternative: either a subgroup is virtually abelian or it contains a rank-two free group \cite{Mc}, \cite{I}.  The latter is a sufficient condition for exponential growth.  In fact, the mapping class group has uniform exponential growth, because its action on homology yields a surjection onto a linear group with uniform exponential growth \cite{AAS}.  Here we prove that its subgroups with exponential growth have uniform exponential growth, and furthermore their minimal growth rates have a lower bound determined by the surface.  This is a consequence of a stronger fact; we state both precisely in the following theorem.

Where $A$ generates the group $G$, let $A$-length denote the length of an element of $G$ in the word metric induced by $A$.

\begin{thm}\label{uueg} Suppose $S$ is a compact orientable surface with mapping class group $\Mod{S}$ and $G < \Mod{S}$ is finitely generated and not virtually abelian.
\begin{itemize}
\item \emph{(Short independent words)} There exists a constant $w$ depending only on $S$ such that, for any finite set $A$ generating $G$, some pair of elements with $A$-length less than $w$ generates a rank two free group.
\item \emph{(Uniform uniform exponential growth)} There exists a positive constant $r$ depending only on $S$ such that $h(G) > r$.
\end{itemize}
\end{thm}

Breuillard pointed out the independent interest of the first statement, which he recently proved in the linear group setting, where $G$ is a non-virtually-solvable subgroup of $GL_n(K)$ and $w$ depends only on $n$ \cite{Br}.  One cannot hope for a bound on $w$ independent of dimension, in the linear group setting, or surface complexity, in the mapping class group setting.  This is because there exist linear groups with arbitrarily small minimal growth rates (\cite{GH}; see Remark 1.4 in \cite{Br}).  Furthermore, these same examples embed into mapping class groups, by a construction suggested by Breuillard and described in the last section.  Thus we have

\begin{thm} \label{slowgroups} Let $w = w(S)$ and $r = r(S)$ be the constants from Theorem \ref{uueg}.  There exists a sequence of surfaces $S_n$ such that $w(S_n) \goesto \infty$ and $r(S_n) \goesto 0$ as $n \goesto \infty$. \end{thm}

One application of the first part of Theorem \ref{uueg} appears as a comment in \cite{Br}.  The existence of short independent words yields a fact about simple random walks: return probability decays exponentially with a rate depending only on $w$ and the initial support of the walk.  In other words, start at the identity and iterate a uniform probability distribution on a symmetric set (i.e., one including elements and their inverses) $A$ to obtain a simple random walk on the group $\grp{A}$.  Let $p^{(n)}$ be the probability that one returns to the identity in $n$ steps.  The quantity $\rho = \lim_{n \to \infty}(p^{(n)})^{\frac{1}{n}}$ is the Kesten spectral radius, introduced in \cite{Ke}, and regarding this one can show the corollary below, whose proof we include in Section \ref{proofs}.

\begin{cor}\label{rwalk} There exists a function $f : \mathbb{N} \to (0,1)$ depending only on $S$ which fulfills the following: for the simple random walk described above where $A \subset \Mod{S}$, $|A| = k$, and $\grp{A}$ is not virtually abelian, one has $\rho \leq f(k)$.
\end{cor}

Our proof of Theorem \ref{uueg} relies on Thurston's classification of mapping classes and Ivanov's description of subgroups of the mapping class group to consider three kinds of generators: pseudo-Anosov maps, Dehn twists, and maps which restrict to pseudo-Anosovs on some subsurface.  Specifically, we obtain a constant $p$ with the property that, for any mapping classes $a$ and $b$ in a certain finite-index subgroup of $\Mod{S}$, if $\grp{a,b}$ contains the free group $F_2$, at least one of the following situations applies:
\begin{itemize}
\item[(a)] If one, say $a$, is pseudo-Anosov, then $\grp{a^p,ba^{p}b^{-1}} \cong F_2$.
\item[(b)] If both $a$ and $b$ are Dehn twists, then $\grp{a^p,b^p} \cong F_2$.
\item[(c)] Otherwise $a$ or $b$ restricts to a pseudo-Anosov on a proper subsurface;\\in this case either $\grp{a^p,b^p}$, $\grp{a^p,b^{p}a^{p}b^{-p}}$, $\grp{b^p,a^{p}b^{p}a^{-p}}$, or $\grp{a^p,ba^{p}b^{-1}} \cong F_2$
\end{itemize}

Case (a) is due to K. Fujiwara \cite{F}.  To establish case (b), we generalize a result of Hamidi-Tehrani \cite{H}.  Case (c) is our primary contribution.

\medskip

\noindent {\bf Acknowledgments.}  Many thanks to Juan Souto for posing and pushing the question; Dick Canary for enduring the author's rambles, re-writes, and stalls; and Chris Leininger for the path to enlightenment.  Most of this work was completed at MSRI's Workshop on Teichm\"{u}ller Theory and Kleinian Groups, and the author is grateful for the generosity and interest of many of its participants, including: Javier Aramayona, Moon Duchin, Koji Fujiwara, Aditi Kar, and Hossein Namazi.  Much appreciation also to Emmanuel Breuillard for his comments and insight.

\section{Preliminaries}

\subsection{Mapping class group elements and subgroups.}

Let $S = S_{g,p}$ be an oriented surface with genus $g$, punctures $p$, and \emph{complexity} $\xi(S) = 3g + p$. Its \emph{mapping class group} $\Mod{S}$ consists of orientation-preserving homeomorphisms of $S$ up to isotopy.  Elements of this group are \emph{mapping classes}.

Let \emph{curve} be shorthand for the isotopy class of a simple closed curve on $S$ that bounds neither a disk nor an annulus.  The \emph{intersection number} $i(\gamma,\gamma')$ of two curves $\gamma, \gamma'$ is the minimum number of points of intersection over pairs of representatives of their isotopy classes; if $i(\gamma,\gamma') = 0$ the curves are \emph{disjoint}, else they \emph{intersect}.  This definition extends to \emph{multicurves}, i.e. nonempty collections of disjoint, distinct curves $\gamma = \{\gamma_i\}_{i=1}^{n}$.

The simplest example of an infinite-order mapping class is a \emph{(right) Dehn twist} about a curve $\gamma$: parametrize an annular neighborhood of $\gamma$ with (orientation-preserving) coordinates  $\{(x,y) \in \R \times [0,1] \}/(x,y)\sim(x+1,y)$, and apply the map defined by $(x,y) \mapsto (x+y,y)$ inside the annulus and the identity everywhere else.  In this note we also use the term to describe any composition of Dehn twists about disjoint curves.  The following lemma, which we need later, may be used to prove that Dehn twists have infinite order.  Let $T_\gamma$ denote the Dehn twist about $\gamma$.

\begin{lmm}[Intersection after Dehn twisting]\label{intersection} Let $\gamma_1,\dots,\gamma_s$ be disjoint curves and let $T = T^{e_1}_{\gamma_1} \cdots T^{e_s}_{\gamma_s}$.  Then for all multicurves $\delta, \delta'$ on $S$,
\[ i(T(\delta),\delta') \geq \sum_j(|e_j|-2)i(\delta,\gamma_j)i(\delta',\gamma_j)-i(\delta,\delta') \] \end{lmm}

\noindent This is Lemma 4.2 in \cite{I}.  If the $e_j$ have the same sign, the $- 2$ may be omitted.

Note that $T_\gamma$ fixes $\gamma$; in general, if a mapping class fixes a multicurve setwise it is called \emph{reducible}.  Thurston characterized the infinite-order irreducible mapping classes as those that stretch and shrink a pair of singular foliations of the surface; these are called \emph{pseudo-Anosov} in analogy with Anosov maps on the torus.  Given a reducible mapping class $\phi$, we want to cut the surface along its fixed multicurve so that $\phi^k$ acts like the identity map or a pseudo-Anosov on the components, where the power $k$ accounts for the possibility that $\phi$ permutes the components.  Thurston's classification of $\Mod{S}$ elements enables such a decomposition, and Birman, Lubotzky, and McCarthy proved its uniqueness \cite{FLP}, \cite{BLM}.  Ivanov generalized the ideas to apply to subgroups instead of individual elements, obtaining the results we paraphrase for the remainder of this subsection \cite{I}.

For our purposes we may restrict attention to a particular finite-index subgroup of $\Mod{S}$, in which the problem of permuting components does not arise.  Let $\Gamma(S)$ be the kernel of the $\Mod{S}$-action on homology with $\Z/3\Z$ coefficients; call the elements of $\Gamma(S)$ \emph{pure} mapping classes.

\begin{thm}[Structure of pure mapping classes]\label{thurston} Associated to a pure mapping class $\phi$ is a unique multicurve $C$ such that \begin{itemize}
\item $\phi$ has a representative homeomorphism $f$ which fixes a representative of each component of $C$,
\item where $R_1,\dots R_n$ are the components of $S$ cut along $C$, $f(R_i) = R_i$, and $f|_{R_i}$ is either pseudo-Anosov or the identity, and
\item $i(\gamma,C) > 0$ implies $\phi(\gamma) \neq \gamma$.
\end{itemize}
\end{thm}
\noindent Call $C$ a \emph{canonical reducing system for $\phi$}, and the $R_i$ \emph{pseudo-Anosov components} or \emph{identity components} depending on $f|_{R_i}$.  If $C$ is empty, $\phi$ is either the identity or pseudo-Anosov.  If $C$ is nonempty but all the $R_i$ are identity components, $\phi$ is a Dehn twist.  Otherwise $\phi$ has a pseudo-Anosov component $R_i \neq S$; in this case we call $\phi$ \emph{relatively pseudo-Anosov}, borrowing terminology from Hamidi-Tehrani \cite{H}.  We also consider the components $R_i$ as subsurfaces of $S$; then the theorem above may be restated as: every pure mapping class is either a Dehn twist, pseudo-Anosov, or relatively pseudo-Anosov, the latter meaning it restricts to a pseudo-Anosov on a proper subsurface.

Notions for $\Mod{S}$ elements extend to $\Mod{S}$ subgroups.  An \emph{irreducible subgroup} is one for which no multicurve is preserved by all elements of the group.  One may associate to a subgroup $H$ a canonical reducing system $C$, and subsurfaces $R_i$ from $S$ cut along $C$, such that one has induced groups $H_i < \Mod{R_i}$ which are each either trivial or irreducible.  For a pure subgroup $H < \Gamma(S)$, the $H_i$ are generated by the restrictions $h|_{R_i}$ of elements of $H$ to the subsurface $R_i$; each $H_i$ is also pure, and furthermore the nontrivial $H_i$ are either infinite cyclic or contain a rank-two free group.  If the latter occurs for some $H_i$, then $H$ itself contains a rank-two free group; otherwise $H$ is abelian---we will refer to this fact as the ``strong Tits alternative'' for pure subgroups.  We state the other relevant consequences in a theorem for later reference:

\begin{thm}\label{ivanov}If $H < \Gamma(S)$ contains a free group, then on some subsurface $R$ of $S$, $H$ induces an irreducible subgroup $H' < \Gamma(R)$ which contains a free group.  If $a,b \in H$ restrict to $a_R,b_R \in H_R$ such that $\grp{a_R,b_R} \cong F_2$, then $\grp{a,b} \cong F_2$.\end{thm}

\subsection{The curve complex and projection of curves to subsurfaces.}

The \emph{curve complex} $\C{S}$ of a surface $S$ gives a way to quantify the action of mapping classes on curves.  For a surface with complexity $\xi(S) > 4$, $\C{S}$ has vertex set $\V{S}$ representing the connected curves on $S$, and edges joining vertices representing disjoint curves.  Higher simplices correspond to n-tuples of pairwise disjoint curves, but we only need the 1-skeleton $\G{S}$. We are also interested in the cases where $\xi(S) = 4$, namely the once-punctured torus $S_{1,0}$ and the four-punctured sphere $S_{0,4}$.  For these, the vertices are still curves on the surface but edges join curves which intersect minimally for any two curves on the surface---i.e., once for $S_{1,0}$ and twice for $S_{0,4}$.  Both $\G{S_{1,0}}$ and $\G{S_{0,4}}$ are the familiar Farey graph.  Give $\G{S}$ the path metric where edges have unit length, and let $\dis{S}{\gamma,\gamma'}$ denote the distance between vertices $\gamma,\gamma'$ in $\G{S}$.  $\Mod{S}$ acts on $\G{S}$ by isometries.  Of all the interesting facts about the action of pseudo-Anosovs on $\C{S}$, we need the following:

\begin{thm}[Minimal translation of pseudo-Anosovs]\label{masurminsky} Suppose $\xi(S) \geq 4$.  There exists $c > 0$ depending only on $S$ such that, for any pseudo-Anosov $h \in \Mod{S}$, any curve $\gamma \in \V{S}$, and any $n \in \Z-\{0\}$, \[\dis{S}{h^n(\gamma),\gamma} \geq c|n| .\] \end{thm}

\noindent This theorem for $\xi(S) > 4$ comes from Masur and Minsky in \cite{MM1}.  For $\G{S_{1,0}}$ and $\G{S_{0,4}}$ one obtains the same result by considering hyperbolic actions on the Farey graph.

We also want to relate $\V{S'}$ to $\V{S}$ where $S'$ is a subsurface of $S$ and $\xi(S),\xi(S') \geq 4$.  We always assume that a subsurface $S' \subset S$ is connected and essential, i.e. no component of $\del S'$ bounds a disk in $S$.  As with curves, subsurfaces are only defined up to isotopy; this poses no problem in light of Theorem \ref{thurston}.  Note that only surfaces with complexity $\xi(S) \geq 4$ can support a pseudo-Anosov.

Here we recall ideas described at length in \cite{MM2}.  Define a projection $\pi_{S'}$ from $\V{S}$ to subsets of $\V{S'}$ as follows: represent $\gamma \in \V{S}$ by a curve with minimal intersection with $\del S'$.  If $\gamma \subset S'$, $\pi_{S'}(\gamma) = \{\gamma\}$.  Otherwise, for each arc $\gamma_i$ of $\gamma \cap S'$, take the boundary of a regular neighborhood of $\gamma_i \cup \del S'$ and exclude the component curves which bound an annulus in $S'$; call the remaining curves $\gamma'_i$.  Then $\pi_{S'}(\gamma) = \cup_{i}\gamma'_i$.  If $\gamma$ misses $S'$, then $\pi_{S'}(\gamma)$ is the empty set; otherwise it is a subset of $\V{S'}$ with diameter at most two, as shown in Lemma 2.2 of \cite{MM2}.  If $\pi_{S'}(\alpha)$ and $\pi_{S'}(\beta)$ are non-empty we may define the \emph{projection distance} $\dis{S'}{\alpha,\beta}$ as the diameter in $\V{S'}$ of $\pi_{S'}(\alpha) \cup \pi_{S'}(\beta)$.

Let us now present a useful lemma about subsurface projection, due to Behrstock.  To keep constants as constructive as possible, we record the relatively elementary proof indicated to me by Chris Leininger.  Call two subsurfaces $A$ and $B$ of $S$ \emph{overlapping} if they are neither disjoint nor nested (equivalently, both $\pi_{A}(\del B)$ and $\pi_{B}(\del A)$ are non-empty).

\begin{lmm} [Behrstock \cite{B}] \label{behrstock}For any pair of overlapping subsurfaces $Y$ and $Z$ of $S$ such that $\xi(Y), \xi(Z) \geq 4$, and any multicurve $x$ with nonempty projection to both,
\[ \dis{Y}{x,\del Z} \geq 10 \sez \dis{Z}{x,\del Y} \leq 4 \]
\end{lmm}

This follows after a few facts.  Suppose $S'$ is a subsurface of $S$ and $\xi(S),\xi(S') \geq 4$.  Let $u_0$ and $v_0$ be curves on $S$ which minimally intersect $S'$ in sets of arcs which include $a_u$ and $a_v$ respectively, and suppose $u$ is a component of the boundary of a neighborhood of $a_u \cup \del S'$, as is $v$ for $a_v \cup \del S'$, so that $u \in \pi_{S'}(u_0)$ and $v \in \pi_{S'}(v_0)$.  Define intersection number of arcs to be minimal over isotopy fixing the boundary setwise but not necessarily pointwise.  One has:

\begin{itemize}
\item Fact 1. If $i(a_u,a_v) = 0$, then $\dis{S'}{u_0,v_0} \leq 4$
\item Fact 2. If $i(u,v) > 0$, then $i(u,v) \geq 2^{(\dis{S'}{u,v} -2)/2}$
\item Fact 3. $i(u,v) \leq 2 + 4 \cdot i(a_u,a_v)$
\end{itemize}

Fact 1 follows from Lemma 2.2 in \cite{MM2}.  A straightforward induction proves Fact 2, which Hempel records as Lemma 2.1 in \cite{He}.  Fact 3 is the observation that essential curves from the regular neighborhoods of $a_u \cup \del S'$ and $a_v \cup \del S'$ intersect at most four times near every intersection of $a_u$ and $a_v$, plus at most two more times near $\del S'$.

\medskip

\noindent \textbf{Proof of Behrstock's Lemma.}  Since $\dis{Y}{x,\del Z} \geq 10 > 2$, the diameter is realized by curves $u \in \pi_{Y}(x), v \in \pi_{Y}(\del Z)$ such that, by Fact 2, $i(u,v) \geq 2^4 = 16$.  By the definition of $\pi_Y$, these $u$ and $v$ come from arcs $a_u \subset x \cap Y$ and $a_v \subset \del Z \cap Y$ respectively.  By Fact 3, $i(a_u,a_v) \geq (16-2)/4 > 3$.  Thus $a_u$ is an arc of $x$ intersected thrice by an arc $a_v$ of $\del Z$, within the subsurface $Y$.  Observe that one of the segments of $a_u$ between points of intersection must lie within $Z$.  This segment is an arc $a_x$ of $x$ disjoint from arcs of $\del Y$ in $Z$.  Fact 1 implies $\dis{Z}{x,\del Y} \leq 4$.\QED

\subsection{Growth bounds and ping-pong}  Two elements $a,b$ of a group $G$ are \emph{independent} if $\grp{a,b}$ is a rank-two free subgroup of $G$.  If one finds independent elements of length less than $d$ in the word metric induced by some generating set $A$, $h(G,A) \geq (\log 3)/d$.  To prove independence, one has the following lemma, a geometric group theory standard.

\begin{lmm}[Ping-pong]\label{pingpong} Suppose $G=\grp{a,b}$ acts on a set $X$, and suppose there exist nonempty disjoint subsets $X_a,X_b \subset X$ such that $a^k(X_b) \subset X_a$ and $b^k(X_a) \subset X_b$ for all $k \in \Z - \{0\}$. Then $G$ is a rank-2 free group.\end{lmm}

\noindent \textbf{Proof.} Any non-empty reduced word of the form $w = a^*b^*\cdots b^*a^*$ (*'s are non-zero integers) cannot be the identity because $w(X_b) \cap X_b \subset X_a \cap X_b = \nada$.  Any other reduced word is conjugate to $w$ in the above form. \QED

\section{Short words, fast growth}\label{proofs}

\noindent Theorem \ref{uueg} follows from:

\begin{lmm}[Main lemma]\label{main} There exists a constant $p$ depending only on $S$ such that, for any pure mapping classes $a,b \in \Mod{S}$, if $\grp{a,b}$ contains the free group $F_2$, then one of the following holds:
\begin{itemize}
\item[(a)] If one, say $a$, is pseudo-Anosov, then $\grp{a^p,ba^{p}b^{-1}} \cong F_2$.
\item[(b)] If both $a$ and $b$ are Dehn twists, then $\grp{a^p,b^p} \cong F_2$.
\item[(c)] Otherwise $a$ or $b$ is relatively pseudo-Anosov; in this case either\\$\grp{a^p,b^p}$, $\grp{a^p,b^{p}a^{p}b^{-p}}$, $\grp{b^p,a^{p}b^{p}a^{-p}}$, or $\grp{a^p,ba^{p}b^{-1}} \cong F_2$
\end{itemize}
\end{lmm}

\noindent \textbf{Proof of Theorem \ref{uueg}.} Suppose $G<\Mod{S}$ has finite generating set $A$.  Consider its finite-index pure subgroup $H = G \cap \Gamma(S)$.  Lemma 3.4 of Shalen and Wagreich in \cite{SW} gives a finite generating set $A'$ for $H$, whose elements have $A$-length at most $2[G:H]-1$.  If $G$ is not virtually abelian, $H$ is not abelian, so one finds a noncommuting pair $a,b \in A'$.  By the ``strong Tits alternative'' for pure mapping classes described at the end of 2.1, $\grp{a,b}$ contains a free group.  Then the main lemma gives two independent elements of $A'$-length bounded above by $3p$, thus $A$-length bounded by $3p \cdot (2[G:H]-1)$.  It is straightforward to check that one may use $$w = 3p \cdot (2[\Mod{S}:\Gamma(S)]-1) \qquad \text{and} \qquad r = (\log 3)/w$$ for Theorem \ref{uueg}.\QED

The main lemma compiles three separate results.  For case (a), we turn to work of Fujiwara regarding group actions which are ``acylindrical,'' a property Bowditch defined and proved for the $\Mod{S}$-action on $\C{S}$ \cite{Bo}.  The statement we need follows from Theorem 3.1(2) in \cite{F}:

\begin{thm}[Pseudo-Anosovs]\label{fujiwara} There exists a power $p_0$ depending only on $S$ with the following property: let $a$ be a pseudo-Anosov mapping class and $b$ any mapping class such that $\grp{a,b}$ contains $F_2$.  Then for any $k > p_0$, $\grp{a^{k},ba^kb^{-1}}$ is a rank-2 free group. \end{thm}

To obtain case (b), we observe that a result of Hamidi-Tehrani in \cite{H}, concerning Dehn twists where all powers have the same sign, extends to the statement we need if one uses Lemma \ref{intersection} to estimate intersection numbers after Dehn twisting, instead of the sharper inequality available in the same-sign case.  We repeat the details here because they mirror our proof for case (c).

\begin{prp}[Dehn twists]\label{ht}  Suppose two Dehn twists $a$ and $b$ do not commute.  Then for any $p \geq 4$, $\grp{a^p,b^p}$ is a free group. \end{prp}

\noindent \textbf{Proof.}  Assume $a = T^{k_1}_{\alpha_1}T^{k_2}_{\alpha_2} \cdots T^{k_n}_{\alpha_n}$ and $b = T^{l_1}_{\beta_1}T^{l_2}_{\beta_2} \cdots T^{l_m}_{\beta_m}$, where $\alpha = \{\alpha_i\}$ and $\beta = \{\beta_i\}$ are multicurves, and the $k_i$ and $l_j$ are nonzero powers.  Since $a$ and $b$ do not commute, $\alpha$ intersects $\beta$; by restricting attention to a subsurface we may further assume that $i(\alpha,\beta_j)$ and $i(\alpha_i,\beta)$ are nonzero for all $i,j$.  Under the additional assumption that $|k_i| \geq 4$ and $|l_j| \geq 4$ for all $i,j$, we prove that $\grp{a,b}$ is free; this is slightly more general than the proposition stated above.

Recall $\V{S}$ represents curves on $S$.  We want to apply the ping-pong lemma.  Define ping-pong sets as follows:
\begin{eqnarray*}
X_a & = & \{\gamma \in \V{S} : i(\gamma,\alpha) < i(\gamma,\beta)\}  \\
X_b & = & \{\gamma \in \V{S} : i(\gamma,\beta)  < i(\gamma,\alpha)\} \\
\end{eqnarray*}
\noindent  It is clear $X_a \cap X_b = \nada$, and our assumptions ensure that $X_a$ and $X_b$ are nonempty.  It remains to check that for all $k \in \Z - \{0\}$, $a^k(X_b) \subset X_a$ and $b^k(X_a) \subset X_b$.  Suppose $\gamma \in X_a$.  Applying Lemma~\ref{intersection},
\begin{eqnarray*}
i(b^k(\gamma),\alpha) &\geq& \sum_j(|l_j|-2)i(\gamma,\beta_j)i(\alpha,\beta_j)- i(\gamma,\alpha) \\
                      & > &   \sum_j(|l_j|-2)i(\gamma,\beta_j)i(\alpha,\beta_j)- i(\gamma,\beta) \\
                      & = &   \sum_j((|l_j|-2)i(\alpha,\beta_j)-1)i(\gamma,\beta_j) \\
                      &\geq& \sum_j((4-2)(1)-1)i(\gamma,\beta_j) \\
                      & = & i(\gamma,\beta) = i(b^k(\gamma),b^k(\beta)) \\
                      & = & i(b^k(\gamma),\beta)
\end{eqnarray*}  Thus $b^k(X_a) \subset X_b$.  By symmetry, $a^k(X_b) \subset X_a$.\QED

Notice that, in the preceding proof, $a$ and $b$ fix multicurves $\alpha$ and $\beta$ respectively, $b$ takes curves which predominantly intersect $\beta$ to curves which predominantly intersect $\alpha$, and $a$ does the opposite, so ping-pong proves $\grp{a,b}$ free.  For case (c), we consider $a$ and $b$ acting as pseudo-Anosovs on proper subsurfaces $A$ and $B$ respectively.  We find a similar set-up for ping-pong: $a$ and $b$ fix the multicurves $\del A$ and $\del B$ respectively, $b$ takes curves which intersect $\del B$ many times in the subsurface $A$ to curves which intersect $\del A$ many times in the subsurface $B$, and $a$ does the opposite.  This idea leads to:

\begin{prp}[Relative pseudo-Anosovs]\label{johanna} There exists a constant $p_1$ depending only on $S$ such that the following holds: suppose $a$ and $b$ are pure mapping classes, $a$ is relatively pseudo-Anosov, and $\grp{a,b}$ is irreducible.  Then either $\grp{a^k,b^k}$, $\grp{a^k,b^{k}a^{k}b^{-k}}$, or $\grp{b^k,a^{k}b^{k}a^{-k}}$ is a rank-two free group for any $k > p_1$. \end{prp}

\noindent \textbf{Proof.}  The first possibility, $\grp{a^k,b^k} \cong F_2$, happens when $a$ and $b$ have overlapping pseudo-Anosov components; we prove this before addressing the general situation.

Let $A$ and $B$ be overlapping pseudo-Anosov components of $a$ and $b$ respectively, which means that $\del A$ projects nontrivially to $B$, as $\del B$ does to $A$.  Let $\Sigma$ be the finite set of topological types of essential subsurfaces of $S$ with $\xi \geq 4$ (i.e. all possible pseudo-Anosov components, including $S$ itself) and for any $S' \in \Sigma$ let $c(S')$ be the constant from Theorem \ref{masurminsky}.  Set

\[ p_1 = \max [\{ 14/c(S') : S' \in \Sigma \} \cup \{2\}]. \]

We use the ping-pong lemma to show that for any $k > p_1$, $\grp{a^k,b^k}$ is a free group.  Define ping-pong sets as follows:
\begin{eqnarray*}
X_a & = & \{\gamma \in \V{S} : \pi_A(\gamma) \neq \nada, \pi_B(\gamma) \neq \nada,
                           \dis{A}{\gamma,\del B} \geq 10 \}  \\
X_b & = & \{\gamma \in \V{S} : \pi_A(\gamma) \neq \nada, \pi_B(\gamma) \neq \nada,
                           \dis{B}{\gamma,\del A} \geq 10 \}
\end{eqnarray*}

Lemma \ref{behrstock} guarantees $X_a \cap X_b = \nada$.  Suppose $\gamma \in X_a$ and $m \in \Z - \nada$.  Applying the same lemma,
\[ \dis{A}{\gamma,\del B} \geq 10 \sez \dis{B}{\gamma,\del A} \leq 4. \]
We chose $p_1$ so that, applying Theorem \ref{masurminsky},
\[ \dis{B}{\gamma,b^{mk}(\gamma)} \geq c(B) \cdot |mk| \geq c(B) \cdot p_1 \geq c(B) \cdot 14/c(B) = 14.  \]
The triangle inequality gives
\[ \dis{B}{b^{mk}(\gamma),\del A} \geq \dis{B}{\gamma,b^{mk}(\gamma)} - \dis{B}{\gamma,\del A} \geq 14 - 4 = 10.\]
Thus $b^{mk}(X_a) \subset X_b$.  To show that $X_b$ is nonempty, apply Theorem \ref{masurminsky} to show $b^{k}(\del A) \in X_b$.  A symmetric argument shows $X_a$ is nonempty and $a^{mk}(X_b) \subset X_a$.  Thus Lemma \ref{pingpong} shows that $\grp{a^k,b^k}$ is free.

Now suppose $a$ and $b$ are as in the theorem.  Choose a pseudo-Anosov component $A$ of $a$.  Irreducibility of $\grp{a,b}$ implies $\del A$ cannot be fixed by $b$, so one of the following cases applies:
\begin{itemize}
\item Case 1: Suppose $b$ is pseudo-Anosov.  Observe that for $k > p_1 > 1/c(S)$, $\dis{S}{\del A, b^{k}(\del A)} > 1$.  This means $b^k(\del A)$ and $\del A$ intersect, which guarantees $A$ and $b^{k}A$ are overlapping pseudo-Anosov components of $a$ and $b^{k}ab^{-k}$.  Thus $\grp{a^{k},b^{k}a^{k}b^{-k}} \cong F_2$.
\item Case 2: Suppose $b$ is reducible and $\del A$ lives in the closure of the union of identity components of $b$.  Then $b$ restricted to this union is a Dehn twist.  Applying Lemma \ref{intersection}, we see that for $k > p_1 > 2$, $i(b^k(\del A),\del A)$ is positive.  Again, $A$ and $b^{k}A$ are overlapping pseudo-Anosov components of $a$ and $b^{k}ab^{-k}$, so $\grp{a^{k},b^{k}a^{k}b^{-k}} \cong F_2$.
\item Case 3: Suppose $b$ is reducible and $\del A$ projects nontrivially to some pseudo-Anosov component $B$ of $b$.  Irreducibility of $\grp{a,b}$ implies $\del B$ cannot be fixed by $a$.
    \begin{itemize}
    \item If $\del B$ lives in the closure of the union of identity components of $a$, we are in Case 2 with the roles of $a$ and $b$ reversed. For any $k>p_1$, $\grp{b^{k},a^{k}b^{k}a^{-k}} \cong F_2$.
    \item Otherwise $\del B$ projects to some pseudo-Anosov component $A'$ of $a$, equal to or disjoint from $A$.  Either $A$ and $B$ overlap or $A$ is nested in $B$, in which case $B$ cannot nest in $A'$, so $B$ and $A'$ overlap.  In any case we obtain overlapping pseudo-Anosov components of $a$ and $b$, so for $k>p_1$, $\grp{a^k,b^k} \cong F_2$. \QED
    \end{itemize}
\end{itemize}

\medskip

\noindent \textbf{Proof of Main Lemma.}  The reduction to cases (a), (b), and (c) follows from Theorem \ref{thurston}.  While cases (a) and (b) correspond to Theorem \ref{fujiwara} and Proposition \ref{ht} respectively, case (c) is not identical to Proposition \ref{johanna}, because $\grp{a,b}$ may be reducible.  Nevertheless, Theorem \ref{ivanov} gives a subsurface $R$ such that $\grp{a|_R,b|_R}$ is irreducible and contains $F_2$.  By construction, the constant $p_1$ of Proposition \ref{johanna} works if we replace $S$ with one of its subsurfaces.  However, the same need not hold for the constant $p_0$ of Theorem \ref{fujiwara}, which we must apply if both $a|_R$ and $b|_R$ are pseudo-Anosov.  To remedy this, let $p_2 = \max \{ p_0(S') : S' \in \Sigma \}$, where, as in the proof of Proposition \ref{johanna}, $\Sigma$ represents subsurfaces of $S$.  Then for the constant in the Main Lemma, we may choose any $p > \max\{4, p_1, p_2\}$.  \QED

\noindent \textbf{Proof of Corollary \ref{rwalk}.}  For a reference on random walks, see \cite{Wo}; also, the author enjoys the exposition in \cite{KV}.  For a symmetric set $A$, let $\rho(A)$ be the spectral radius for the corresponding simple random walk on $\grp{A}$, as described in the introduction.  Part of Theorem 10.3 in \cite{Wo}, specialized to simple random walks, gives the following equivalence: $\rho(A)$ is strictly less than 1 if and only if there exists $\kappa(A) > 0$ such that, for all finitely supported functions $f:\grp{A} \to \mathbb{C}$,
\[\sum_{x \in \grp{A}} f(x)^2 \leq \kappa(A)\sum_{x \in \grp{A}, g \in A}(f(x) - f(xg))^2/2.\]
One may think of the sum on the right as the sum of squares of differences over edges in the Cayley graph determined by $A$.  For the forward implication, one may choose $\kappa \geq (1-\rho)^{-1}$; backwards, one knows $\rho \leq 1-\kappa^{-1}$ (note this corrects a typo in the proof in \cite{Wo}).

Suppose $A$ is as in the hypothesis: $A \subset \Mod{S}$, $|A| < k$, and $\grp{A}$ is not virtually abelian.  Theorem \ref{uueg} gives two elements $u,v \in \grp{A}$ with $A$-length less than $w$ such that $u$ and $v$ freely generate $F_2$.  Let $T$ be the symmetric set $\{u,v,u^{-1},v^{-1}\}$.  Then $\rho(T) = \sqrt{3}/2$, since this is the simple random walk on standard generators of $F_2$, already calculated in \cite{Ke}.  From $\rho(T)$ one obtains $\kappa(T)$, and, because cosets of $\grp{T}$ are for this purpose indistinguishable from $\grp{T}$,
\[\sum_{x \in \grp{A}}f(x)^2 \leq \kappa(T)\sum_{x \in \grp{A}, t \in T}(f(x)-f(xt))^2/2\] for any finitely supported function $f:\grp{A} \to \mathbb{C}$.  Furthermore, since the generators in $T$ are of the form $g_{1}\cdots g_{d}$, $g_i \in A$, $d \leq w$, the triangle inequality implies that for any $x \in \grp{A}, t \in T$, there exists $y \in \grp{A}, g \in A$ such that \[(f(x)-f(xt))^2 \leq w^2(f(y)-f(yg))^2.\]  In the Cayley graph, the pair $(y,yg)$ corresponds to an edge of the path connecting $x$ to $xt$.  It is straightforward to count that for a given $y \in \grp{A}$ and $g \in A$, the difference $f(y)-f(yg)$ need appear at most $w(k-1)^{w-1}$ times in the upper bound for some $f(x)-f(xt)$, $x \in \grp{A}, t \in T$.  Therefore one has \[\sum_{x \in \grp{A}}f(x)^2 \leq \kappa(T)w^3(k-1)^{w-1}\sum_{x \in \grp{A}, g \in A}(f(x)-f(xg))^2/2.\]  Thus \[\rho(A) \leq 1 - {\frac{1 - \sqrt{3}/2}{w^3(k-1)^{w-1}}}.\]  The right hand side defines $f(k)$.\QED

\section{Slow-growing groups}

Here we prove Theorem \ref{slowgroups}, that the constants in Theorem \ref{uueg} necessarily depend on the choice of surface.  It suffices to find subgroups of mapping class groups with arbitrarily small minimal growth rates.  In \cite{GH}, Grigorchuk and de la Harpe construct a sequence of groups $G_n$ with the property that $h(G_n) \goesto 0$ as $n \goesto \infty$.  Following a construction suggested by Emmanuel Breuillard, one may represent these groups in mapping class groups.  In fact, the method shows that finite extensions of mapping class subgroups are themselves mapping class subgroups.  Let $B$ be the non-empty boundary of a surface $S$, $\Homeo{S,B}$ denote the group of orientation-preserving homeomorphisms of $S$ fixing $B$ pointwise, and $\Mod{S,B}$ be $\Homeo{S,B}$ modulo equivalence by isotopy keeping $B$ fixed throughout.

\begin{prp}\label{extend}If a finitely generated group $G$ has a finite index subgroup $H$ isomorphic to a subgroup of $\Mod{S,B}$, then $G$ embeds in the mapping class group of a connected surface.
\end{prp}

\noindent \textbf{Proof.}  Throughout, we consider $H$ as a subgroup of both $\Mod{S,B}$ and $G$.  By passing to a finite index subgroup, we may assume $H$ is normal in $G$.  Ultimately we embed $G$ in the mapping class group of a surface $S^*$ modeled on a Cayley graph of $G/H$ with a copy of $S$ for each vertex, and annuli for edges.

Before constructing $S^*$, let us address two technical issues.  First: we need a group of homeomorphisms associated to $H$, which will let us define $S^*$ as a quotient.  Let $F_A$ be the free group on $A$, where $A = \{g_1, g_2, ... g_m\}$ is a finite generating set for $G$.  Let $\bar{H}$ be the pre-image of $H$ under the canonical map $r: F_A \to G$; observe $F_A/\bar{H}$ and $G/H$ are isomorphic finite groups.  As a subgroup of a free group, $\bar{H}$ is also a free group, and so the map $r: \bar{H} \to H$ lifts to $\tilde{r}: \bar{H} \to \Homeo{S,B}$.  Second: without loss of generality we assume that (1) $B$ has at least $2m$ components, and (2) only the identity in $H$ fixes all curves on $S$ (i.e., no element of $H$ is a Dehn twist about components of $B$).  We can ensure this by replacing $S$ with a larger surface whose mapping class group also contains $H$ \cite{PR}.

Now we build $S^*$.  To each generator $g_i$ associate a pair of boundary components $B_i^+$ and $B_i^-$, and fix orientation-reversing homeomorphisms $b_i: B_i^+ \to B_i^-$.  Let $S'$ be the quotient of the disjoint union $F_A \times S$ by the equivalence $(k,x)~\sim~(kh^{-1},\tilde{r}(h)(x))$ for all $h$ in $\bar{H}$.  One obtains a copy of $S$ for each coset $g\bar{H}$, corresponding to each vertex of the Cayley graph of $F_A/\bar{H}$.  Let $B'$ be the image of $F_A \times B$ in $S'$.  For a connected surface, obtain $S^*$ by identifying components of $B'$ as specified by the edges of the Cayley graph: $[(k,x)] \sim [(kg_i, b_i(x))]$ for all $g_i$ and $x \in B_i^+$.

The action of $F_A$ on $F_A \times S$ by $g \cdot (k,x) = (gk,x)$ descends to an action on $S^*$ that gives a homomorphism $p: F_A \to \Homeo{S^*}$.  Let $q$ be the quotient $q:\Homeo{S^*}\to\Mod{S^*}$.  Define a homomorphism $\phi: G \to \Mod{S^*}$ by $\phi(g) = qpr^{-1}(g)$.  Although $r$ is not injective, $\phi$ is well-defined if, for any two words $w$ and $v$ representing the same element in $G$, $qp(wv^{-1})$ is the identity $e$.  Notice $wv^{-1}$ represents the identity, so it is in $\bar{H}$, thus
$$p(wv^{-1})\cdot[(k,x)]=[(wv^{-1}k,x)]=[(k,\tilde{r}(k^{-1}wv^{-1}k)(x))]$$
Since $r$ maps $k^{-1}wv^{-1}k$ to the identity, the lift $\tilde{r}$ maps it to a homeomorphism isotopic to the identity on $S$, thus on components of $S'$, and thus on $S^*$ (because all homeomorphisms and isotopies fix $B'$, the boundary identifications are no obstruction).  Hence $p(wv^{-1})$ is isotopic to the identity, so $qp(wv^{-1}) = e$.

Finally, we show that $\phi$ is injective.  Suppose $\phi(g) = e$ and choose $w \in r^{-1}(g)$.  Because $qp(w) = e$, $p(w)$ is isotopic to the identity.  In particular, $p(w)$ fixes the subsurface of $S^*$ corresponding to the component $[\{e\} \times S]$ of $S'$.  Therefore $w$ represents an element of $H$.  Furthermore, $p(w)$ fixes the isotopy classes of all curves on that subsurface.  Then by condition (2) on $H$ described above, $w$ represents the identity, so $g = r(w) = e$.\QED

\noindent \textbf{Proof of Theorem \ref{slowgroups}.}  As remarked above, we only need to show that the slow-growing groups $G_n$ are subgroups of mapping class groups.  Proposition \ref{extend} applies because each $G_n$ has a finite-index subgroup $H_n$ isomorphic to the direct product of $2^n$ finitely generated free groups \cite{GH}.  The $H_n$ may be realized as mapping class subgroups using pairs of Dehn twists on disjoint subsurfaces of $S_{2^n}$ minus $2^{n+1}$ disks: mapping classes supported on disjoint subsurfaces commute, and a pair of properly chosen Dehn twists generates the rank two free group, which contains all other finitely generated free groups.\QED

\begin{rmk}Recall that the constant $w(S)$ of Theorem \ref{uueg} grows linearly with both the constant $p(S)$ from Lemma \ref{main} and with $[\Mod{S}:\Gamma(S)]$, where $\Gamma(S)$ are the mapping classes acting trivially on homology with $\Z/3\Z$ coefficients.  The argument above takes advantage of only the latter factor.  This is because, where $G$, $H$, and $S^*$ are as in the proof, one may check $G \cap \Gamma(S^*) < H$.  For $G_n$, the corresponding subgroups $H_n$ have large growth; the minimal growth rates of $G_n$ still approach zero because $[G_n:H_n]$ grows large.  By construction, the constant $p(S)$ from Lemma \ref{main} does not decrease with genus, but one might ask whether it necessarily increases.  The answer is yes if there exists a sequence of surfaces $S_n$ and groups $\Gamma_n < \Gamma(S_n)$ such that $h(\Gamma_n) \goesto 0$ as $n \goesto \infty$.\end{rmk}

\end{document}